\documentclass[10pt,twoside,a4paper]{article}
\usepackage[utf8]{inputenc}
\usepackage[T1]{fontenc}
\usepackage{amsmath}
\usepackage{amsfonts}
\usepackage{amssymb}
\usepackage{amsthm}
\usepackage[pdftex]{graphicx}
\usepackage[bookmarksopen=false,breaklinks=true,%
      backref=page,pagebackref=true,plainpages=false,%
      hyperindex=true,pdfstartview=FitH,colorlinks=false,%
      pdfpagelabels=true,linkcolor=blue,%
      citecolor=red,urlcolor=red,hypertexnames=false%
      ]%
   {hyperref}

\usepackage{mathrsfs}
\usepackage{mathtools}
\usepackage{stmaryrd}
\usepackage[all]{xy}
\usepackage{makeidx}
\usepackage[english]{babel}
\usepackage{xspace}

\usepackage{framed}
\newenvironment{cadre}{
  \framed\selectlanguage{english}}{\endframed}

\renewcommand\paragraph[1]{

\rdb\addcontentsline{toc}{subsubsection}{#1} \medskip \noindent $\bullet$ \textbf{#1}}

\newcommand\gui[1]{``{#1}''}

\newcommand \paref[1] {page~\pageref{#1}}

\newcommand \ssi {if and only if\xspace}

\newcommand \Amo {$A$-module\xspace}

\newcommand \kevs {$k$-vector spaces\xspace}

\newcommand \klg {$k$-\alg}
\newcommand \klgs {$k$-\algs}

\newcommand \agq {algebraic\xspace}
\newcommand \agqt {algebraically\xspace}

\newcommand \alg {algebra\xspace}
\newcommand \algs {algebras\xspace}

\newcommand \algo{algorithm\xspace}
\newcommand \algos{algorithms\xspace}

\newcommand \algq{algorithmic\xspace}

\newcommand \cdi{discrete field\xspace}

\newcommand \coh {coherent\xspace}
\newcommand \cohc {coherence\xspace}

\newcommand \demo {proof\xspace}

\newcommand \dfn{definition\xspace}  
\newcommand \dfns{definitions\xspace}

\newcommand \elts{elements\xspace}

\newcommand \eqv  {equivalent\xspace}

\newcommand \fdi{strongly discrete\xspace} 

\newcommand \id {ideal\xspace}
\newcommand \ids {ideals\xspace}

\newcommand \idep {prime \id}

\newcommand \idps {principal \ids}

\newcommand \itf {\tf \id}
\newcommand \itfs {\tf \ids}

\newcommand \mpfs {\pf modules\xspace}

\newcommand \mtf {\tf module\xspace}
\newcommand \mtfs {\tf modules\xspace}

\newcommand \ncrt{necessarily\xspace}

\newcommand \noet {Noetherianity\xspace}

\newcommand \noe {Noetherian\xspace}

\newcommand \pb{problem\xspace}  
\newcommand \pbs{problems\xspace}

\newcommand \pf {finitely presented\xspace}

\newcommand \pol {polynomial\xspace}
\newcommand \pols {polynomials\xspace}

\newcommand \prt {property\xspace}
\newcommand \prts {properties\xspace}

\newcommand \tf {finitely generated\xspace}

\newcommand \Tho {Theorem\xspace}
\newcommand \tho {theorem\xspace}
\newcommand \thos {theorems\xspace}

\newcommand \cov {constructive\xspace}
\newcommand \cof {\cov}

\newcommand \coma {\cov \maths}
\newcommand \clama {classical \maths}

\renewcommand \cot {constructively\xspace}

\newcommand \maths {mathematics\xspace}

\newcommand {\junk}[1]{}

\renewcommand \leq{\leqslant}

\renewcommand \geq{\geqslant}

\newcommand\tsbf[1]{\textsf{\textbf{\textup{#1}}}}

\newcommand\Tsbf[1]{\hyperref[Ax#1]{\tsbf{#1}}}

\newcommand\sa[1]{\hbox{\large\usefont{T1}{pzc}{m}{it}#1}\,}

\newcommand\Sa[1]{\hyperref[theorie#1]{\sa{#1}}}

\makeatletter 
\def\revddots{\mathinner{\mkern1mu\raise\p@ 
\vbox{\kern7\p@\hbox{.}}\mkern2mu 
\raise4\p@\hbox{.}\mkern2mu\raise7\p@\hbox{.}\mkern1mu}} 
\makeatother

\newcommand \gA {\mathbf{A}}

\newcommand \gR {\mathbf{R}}
\newcommand \gQ {\mathbf{Q}}

\newcommand \gZ {\mathbf{Z}}

\newdimen\xyrowsp
\newcommand{\SCO}[6]{
\xymatrix @R = \xyrowsp {
                                  &1 \ar@{-}[dl] \ar@{-}[dr] \\
#3 \ar@{-}[ddr]                   &   & #6 \ar@{-}[ddl] \\
                                  &\bullet\ar@{-}[d] \\
                                  &\bullet   \\
#2 \ar@{-}[ddr] \ar@{-}[uur]      &   & #5 \ar@{-}[ddl] \ar@{-}[uul] \\
                                  &\bullet \ar@{-}[d] \\
                                  &\bullet  \\
#1 \ar@{-}[uur]                   &   & #4 \ar@{-}[uul] \\
                                  & 0 \ar@{-}[ul] \ar@{-}[ur] \\
}
}

\newcommand \rad {\MA{\mathrm{rad}}}

\newcommand\MA[1]{\mathop{#1}\nolimits}

\newcommand \Sac {\MA{\mathsf{Fsac}}}

\newcommand \Sace {\MA{\mathsf{Fsace}}}

\newcommand \Xn {X_1,\ldots,X_n}

\begin{document}

\title{Algebra in Bishop's style: \\ some major features of the book\\
\gui{A Course in Constructive Algebra}
\\ by Mines, Richman,
and Ruitenburg  }
\author{Henri Lombardi}
\maketitle

\begin{abstract} 
The book  \gui{A Course in Constructive Algebra} (1988) 
 shows the way of understanding classical basic algebra in a constructive style similar to Bishop's Constructive Mathematics.
 Classical theorems are revisited, with a new flavour, and become much more precise. We are often  surprised to find proofs that are simpler and more elegant than the usual ones. In fact, when one cannot use magic tools as the law of excluded middle, it is necessary to understand what is the true content of a classical proof. Also, usual shortcuts allowed in classical proofs introduce sometimes useless detours. 
In order to understand clearly a \pb, prescience may be a handicap.
\end{abstract}

\tableofcontents

\section*{Introduction}  
The book  \gui{A Course in Constructive Algebra} (1988) 
 shows the way of understanding classical basic algebra in a constructive style similar to Bishop's Constructive Mathematics.
 Classical theorems are revisited, with a new flavour, and become much more precise. We are often  surprised to find proofs that are simpler and more elegant than the usual ones. In fact, when one cannot use magic tools as the law of excluded middle (LEM), it is necessary to understand what is the true content of a classical proof. Also, usual shortcuts allowed in classical proofs introduce sometimes useless detours. 
In order to understand clearly a \pb, prescience may be a handicap.

\section{The reception of the book}

The reception of the book in France 
is even more confidential than that of Bishop's book \cite{B67}. 
I have hardly ever met a French mathematician who has but heard of the existence of the book.

The Computer Algebra community could be expected to be a little more up-to-date since all \thos in [CCA] have a computational content, and could, at least in principle, be implemented in the usual Computer Algebra softwares.

Some years ago I have submitted an article of \cov \alg to the section \gui{Computer Algebra}
of the Journal of Algebra, section whose recommendations to the authors explicitly indicate the interest of the journal for \coma.
What was my surprise when the referee asked me to explain what was the precise meaning of  \gui{or}  in \coma, because he was confused and did not understand some arguments. The article was finally rejected in this section of the Journal of Algebra, apparently because of the impossibility of finding a competent referee.


\medskip Nevertheless I have recently discovered the following article
by Sebastian Posur,
\emph{A constructive approach to Freyd categories}.
\url{https://arxiv.org/abs/1712.03492}

\noindent Here is an excerpt from section 2, ``Constructive category theory''. This article seems to me to be a salutary and expected turning point.

\begin{cadre} \label{quote0}
To present our algorithmic approach to Freyd categories, we chose the language of constructive mathematics (see, e.g., [MRR88]). We did that for the following reasons: the language of constructive mathematics
\begin{enumerate}
\item reveals the algorithmic content of the theory of Freyd categories,

\item is perfectly suited for describing generic algorithms, i.e., constructions not depending on particular choices of data structures,
\item allows us to express our algorithmic ideas without choosing some particular model of computation (like Turing machines)
\item encompasses classical mathematics, i.e., all results stated in constructive mathematics are also valid classically,
\item does not differ very much from the classical language in our particular setup.
\end{enumerate}

In constructive mathematics the notions of data types and algorithms (or operations) are taken as primitives and every property must have an algorithmic interpretation. For example given an additive category $\gA$ we interpret the property

\centerline{$\gA$ has kernels} 

\noindent as follows: we have algorithms that compute for given
\begin{itemize}
\item [$\bullet$] $A, B\in \mathrm{Obj}_\gA$, $\alpha \in \mathrm{Hom}_\gA(A,B)$, an object 
$\mathrm{ker}(\alpha)\in \mathrm{Obj}_\gA$ and a morphism 
$$
 \mathrm{KernelEmbedding}(\alpha) \in \mathrm{Hom}_\gA(\mathrm{ker}(\alpha), A)
$$
for which $\mathrm{KernelEmbedding}(\alpha) \cdot \alpha = 0$,
\item [$\bullet$] $A,B,T\in \mathrm{Obj}_\gA$, $\alpha\in \mathrm{Hom}_\gA(A,B)$, $\tau \in \mathrm{Hom}_\gA(T,A)$ such that $\tau \cdot\alpha = 0$ a morphism
$u\in \mathrm{Hom}_\gA(T,\mathrm{ker}(\alpha))$ such that
$$
u \cdot \mathrm{KernelEmbedding}(\alpha) = \tau ,
$$
where $u$ is uniquely determined (up to $=$) by this property. 
\end{itemize}

\smallskip Another important example is given by \emph{decidable equality}, where we interpret the property that for all objects $A, B\in \gA$, we have
$$\forall \alpha, \beta\in \mathrm{Hom}_\gA (A, B)  : (\alpha = \beta) \vee (\alpha \neq  \beta) $$
as follows: we are given an algorithm that decides or disproves equality of a given pair of morphisms\dots 

\smallskip On the other hand, we allow ourselves to work classically whenever we interpret Freyd categories in terms of finitely presented functors. The reason for this is pragmatic: we want to demonstrate the usefulness of having Freyd categories computationally available, and we believe that this can be done by interpreting Freyd categories in terms of other categories that classical mathematicians care about.
\end{cadre}

\section{Revisiting Bishop's set theory}

The authors of [CCA] introduce a philosophy of \maths that differs slightly from that of \cite[Bishop, 1967]{B67}. 
This point of view is probably expressed more directly in the papers \cite{RicConfession,Ric1996} and in the book \cite{Br1987}.

First of all, as in Bishop, the point of view is not that of  formalized mathematics, but of mathematics open to unpredictable developments, and for which the only criterion of truth is the conviction given by a proof.

The mathematical universe is thus not preexisting, it is on the contrary a properly human construction for the use of the human community.

A novelty is the following. The general point of view is to consider that all mathematics, classical as well as constructive, deal with the same ideal objects. The unique difference is in the tools used for the investigation of this universe. Constructive \maths are more general than classical \maths since they use neither LEM nor Choice. Exactly as the theory of groups is more general than the theory of abelian groups, since commutativity is not assumed. 

Let us  quote  a passage.

\begin{cadre} \label{quote1}
Our notion of what constitutes a \textbf{set}%
\index{Set} is a rather liberal one.

\smallskip \noindent \textbf{[I.]2.1 Definition.}
A set $S$ is defined when we describe how to construct its members from
objects that have been, or could have been, constructed prior to $S$, and
describe what it means for two members of $S$ to be equal.

\smallskip Following Bishop we regard the \textbf{equality relation}%
\index{Equality relation} on a set as conventional: something to be
determined when the set is defined, subject only to the requirement that it
be an equivalence relation.

\dotfill\medskip 

A unary relation $P$ on $S$ defines a \textbf{subset}%
\index{Subset} $A=\{x\in S:P(x)\}$ of $S:$ an element of $A$ is an element
of $S$ that satisfies $P$, and two elements of $A$ are equal if and only if
they are equal as elements of $S$. If $A$ and $B$ are subsets of $S$, and if
every element of $A$ is an element of $B$, then we say that $A$ is \textbf{%
contained} in $B$, and write $A\subseteq B$. Two subsets $A$ and $B$ of a
set $S$ are \textbf{equal} if $A\subseteq B$ and $B\subseteq A$; this is
clearly an equivalence relation on subsets of $S$.  \\
We have described how to construct a subset of $S$, and what it means for two subsets of $S$ to be equal. Thus we have defined the set of all subsets, or the \textbf{power set}, of $S$.
\end{cadre}

\smallskip
This is rather surprising for a follower of Bishop. The authors of [CCA] think that the notion of  ``a unary relation defined on a given set'' is so clear that we may consider a well-defined set of all these unary relations. 
In other words, we know how to construct these unary relations, in a similar way as for example we know how to construct a nonnegative integer, or a real number, or a real function. But this seems problematic since nobody thinks that it is possible to have a universal language for mathematics allowing us to codify these relations.
In particular, if the set  $\Omega$ of subsets of the singleton 
$\{0\}$ exists, this means that truth values form a set rather than a class. This seems to say that we know a priori all the truth values that may appear in the future development of \maths.

In fact, it seems that each time a \gui{set of all subsets of \dots} 
is used in the book, this happens in a context where only a well defined set of subsets (in the usual, Bishop, meaning) is necessary. So the set of  \emph{all} subsets is not really needed. Or sometimes the quantification over this \emph{set} is not needed.\footnote{The most important exception is in the definition of well-founded sets and ordinals (see below \paref{ordinals}).}

For example let us see the following \tho, whose proof is incredibly simple and elegant.\footnote{This \tho is not found in classical textbooks. Bourbaki (Algebra, Chapter VII, paragraph 4, section 1), perhaps the best text for this \pb, gives the \tho only for the case  $m=n$, $I_1\neq R$ and $J_1\neq R$. And the \demo is less beautiful than in [CCA].}

\begin{cadre} \label{quote2}
The decomposition in \Tho [V.]2.3 is essentially unique over
an arbitrary commutative ring.

\smallskip \noindent \textbf{[V.]2.4 Theorem.} {\it Let $R$ be a commutative ring, $m\leq n$ positive
integers, and $I_{1}\supseteq I_{2}\supseteq \cdots \supseteq I_{m}$ and $%
J_{1}\supseteq J_{2}\supseteq \cdots \supseteq J_{n}$ ideals of $R$. Suppose 
$M$ is an $R$-module that is isomorphic to $\Sigma _{i=1}^{m}R/I_{i}$ and to 
$\Sigma _{j=1}^{n}R/J_{j}$. Then
\begin{enumerate}
\item[(a)]  $J_{1}=J_{2}=\cdots =J_{n-m}=R$.

\item[(b)]  $I_{i}=J_{n-m+i}$ for $i=1,\ldots ,m$.
\end{enumerate}}
\end{cadre}

Here there is no hypothesis on the \ids $I_i$ and $J_j$. If you would want to formalize completely the discourse, you need the quantification over all ideals of $R$, but you don't really need this complete formalization. 
Similarly, we do not need to quantify over the  class of all commutative rings when we write: \gui{Let $ R $ be a commutative ring}. See \cite[Dependent sums and dependent products in Bishop’s set theory]{petrakis18} for a formal system using class quantification.

\smallskip Note however the following passage which deals with the category of sets, and where the  set $\Omega$ of all subsets of $\{0\}$ plays a crucial role. Note also that the nice  \Tho I.4.1  seems to be mainly aesthetic, without more concrete applications, within the framework of the theory of the categories. 

\begin{cadre} \label{quote3}
[...] The categorical property corresponding to a function $f$ being
one-to-one is that if $g$ and $h$ are maps from any set $C$ to $A$, and $fg=fh$, then $g=h$; that is, $f$ is \textbf{left cancellable}. It is routine
to show that $f$ is one-to-one if and only if it is left cancellable.

A map $f$ from $A$ to $B$ is onto if for each $b$ in $B$ there exists $a$ in 
$A$ such that $f(a)=b$. The corresponding categorical property is that $f$
be \textbf{right cancellable}, that is, if $g$ and $h$ are maps from $B$ to
any set $C$, and $gf=hf$, then $g=h$. The proof that a function $f$ is right
cancellable if and only if it is onto is less routine than the proof of the
corresponding result for left cancellable maps.

\medskip \noindent \textbf{[I.]4.1 Theorem.} \emph{A function is right cancellable in the category of sets if and only if it is onto.}
\begin{proof}
Suppose $f:A\rightarrow B$ is onto and $gf=hf$. If $b\in B$, then there
exists $a$ in $A$ such that $f(a)=b$. Thus $g(b)=g(f(a))=h(f(a))=h(b)$, so $%
g=h$. Conversely suppose $f:A\rightarrow B$ is right cancellable, and let $%
\Omega $ be the set of all subsets of $\{0\}$. Define $g:B\rightarrow \Omega 
$ by $g(b)=\{0\}$ for all $b$, and define $h:B\rightarrow \Omega $\ by 
\[
h(b)=\{x\in \{0\}:b=f(a)%
\text{ for some }a\}.
\]
Thus $h(b)$ is the subset of $\{0\}$ such that $0\in h(b)$ if and only if
there exists $a$ such that $b=f(a)$. Clearly $gf=hf$ is the map that takes
every element of $A$ to the subset $\{0\}$. So $g=h$, whence $0\in h(b)$,
which means that $b=f(a)$ for some $a$.  
\end{proof}
\end{cadre}
 
In fact, an original feature of [CCA] is the consideration of a notion of category as a fully-fledged mathematical object and not as a simple \gui{manière de parler}:

\begin{cadre} \label{quote4}
We deal with two sorts of collections of mathematical objects: sets and
categories. 

\dotfill\medskip 

Given two groups, or sets, on the
other hand, it is generally incorrect to ask if they are equal; the proper
question is whether or not they are \textit{isomorphic}, or, more generally,
what are the homomorphisms between them.

A \textbf{category,%
\index{Category}} like a set, is a collection of objects. An equality
relation on a set constructs, given any two objects $a$ and $b$ in the set,
a \textit{proposition} `$a=b$'. To specify a category $\mathcal{C}$, we must
show how to construct, given any two objects $A$ and $B$ in $\mathcal{C}$, a 
\textit{set} $\mathcal{C}(A,B)$. 
\end{cadre}

A primary interest of categories is to  generalize the notion of a family of objects (indexed by a set). For the category of sets, Bishop \cite{B67}  considers only families of subsets of a given set. But in usual mathematical practice, and particularly in algebra, we sometimes need a more general notion, which corresponds to the notion of dependent types in the constructive theory of types.

\begin{cadre} \label{quote5}
Using the notion of a functor, we can extend our definition of a family of
elements of a set to a family of objects in a category $\mathcal{C}$. Let $I$
be a set. A \textbf{family}%
\index{Family} $A$ \textbf{of objects of $\mathcal{C}$ indexed by} $I$ is a
functor from $I$, viewed as a category, to the category $\mathcal{C}$. We
often denote such a family by $\{A_{i}\}_{i\in I}$. If $i=j$, then the map
from $A_{i}$ to $A_{j}$ is denoted by $A_{j}^{i}$, and is an isomorphism.
\end{cadre}

With these tools, it is possible to construct important objects in today's mathematics, as
\begin{itemize}
\item limits and colimits (e.g. products and coproducts) in some categories,
\item some  \agq structures freely generated by general sets (not \ncrt discrete),
\item many operations on ordinals (see the \dfn of ordinals in  [CCA] below).  
\end{itemize}

\smallskip For example, one proves that a module freely generated by a set $ S $ is flat; but it is not \ncrt projective (Exercise IV.4.9). The classical \tho saying that every module is a quotient of a free module remains valid; the effective consequence is not that the module is a quotient of a projective module, but rather a quotient of a flat module. Thus, by forcing the sets to be discrete (by the aid of LEM),  \clama oversimplify the notion of a free module and lead to conclusions impossible to satisfy algorithmically. 

\smallskip A natural notion of ordinal\footnote{This notion is different from the ones given by Brouwer or Martin-Löf. See also \cite [A constructive theory of ordinals] {CLN2017}.} is also introduced in chapter I of [CCA], and it is used in classification  \pbs of abelian groups (in chapter XI).\label{ordinals}

Note that the definition below of a well-founded set uses the quantification over all subsets of $W$.

\begin{cadre} \label{quote6}
Let $W$ be a set with a relation $a<b$. A subset $S$ of $W$ is said to be 
\textbf{hereditary}  if $w\in S$ whenever $w^{\prime }\in S$ for each 
$w^{\prime }<w$. The set $W$ (or the relation $a<b)$ is \textbf{well founded} if each hereditary subset of $W$ equals $W$. A discrete
partially ordered set is well founded if the relation $a<b$ (that is, $a\leq
b$ and $a\neq b)$ on it is well-founded. An \textbf{ordinal}, or a \textbf{well-ordered set}, is a discrete, linearly ordered, well-founded set.

Well-founded sets provide the environment for arguments by induction. 

\dotfill\medskip 

If $\lambda $ and $\mu $ are ordinals, then an \textbf{injection}%
\index{Injection of ordinals} of $\lambda $ into $\mu $ is a function $\rho $
from $\lambda $ to $\mu $ such that if $a<b$ then $\rho a<\rho b$, and if $%
c<\rho b$, then there is $a\in \lambda $ such that $\rho a=c$. We shall show
that there is at most one injection from $\lambda $ to $\mu .$

\smallskip \noindent \textbf{[I.]6.5 Theorem.} \emph{If $\lambda $ and $\mu $ are ordinals, and $\rho $ and $\sigma$ are injections of $\lambda $ into $\mu $, then $\rho =\sigma .$}

\dotfill\medskip 

\smallskip If there is an injection from the ordinal $\lambda $ to the ordinal $\mu $
we write $\lambda \leq \mu $. Clearly compositions of injections are
injections, so this relation is transitive. By [Theorem]  6.5
it follows that if $\lambda \leq \mu $ and $\mu \leq \lambda $, then $%
\lambda $ and $\mu $ are isomorphic, that is, there is an invertible order
preserving function from $\lambda $ to $\mu $. It is natural to say that two
isomorphic ordinals are \textbf{equal}. 
\end{cadre}
 
We are here in a framework close to the constructive theory of dependent types, where all types are created via inductive \dfns.
 
\section{The corpus of classical abstract algebra treated in the book}

Basic classical algebra is fairly widely covered by the various chapters of [CCA]. Perhaps the best is to recall the table of contents of the book.

\medskip 
\noindent Chapter  I.  Sets.

\noindent 1.  Constructive vs.\ classical mathematics. 2.  Sets, subsets and functions. 3.~Cho\-ice. 4.  Categories. 5.  Partially ordered sets and lattices. 6.  Well-founded sets and ordinals. 7.  Notes.

 \smallskip\noindent  
\noindent Chapter  II.  Basic algebra.

\noindent 1.  Groups. 2.  Rings and fields. 3.  Real numbers. 4.  Modules. 5.  Polynomial rings. 6.  Matrices and vector spaces. 7.  Determinants. 8.  Symmetric polynomials. 9.~Notes.

 \smallskip\noindent  
\noindent Chapter  III.  Rings and modules.

\noindent 1.  Quasi-regular elements and the Jacobson radical. 2.  Coherent and Noetherian modules. 3.  Localization. 4.  Tensor products. 5.  Flat modules. 6.  Local rings. 7.  Commutative local rings. 8.  Notes.

 \smallskip\noindent  
\noindent Chapter  IV.  Divisibility in discrete domains.

\noindent 1.  Divisibility in cancellation monoids. 2.  UFD's and B\'ezout domains. 3.~Dedekind-Hasse rings and Euclidean domains. 4.  Polynomial rings. 5.  Notes. 

 \smallskip\noindent  
\noindent Chapter  V.  Principal ideal domains. 1.  Diagonalizing matrices. 2.  Finitely presented modules. 3.  Torsion modules, $p$-components, elementary divisors. 4.  Linear transformations. 5.  Notes.

 \smallskip\noindent  
\noindent Chapter  VI.  Field theory.

\noindent 1.  Integral extensions and impotent rings. 2.  Algebraic independence and transcendence bases. 3.  Splitting fields and algebraic closures. 4.  Separability and diagonalizability. 5.  Primitive elements. 6.  Separability and characteristic $p$. 7.~Perfect fields. 8.  Galois theory. 9.  Notes.

 \smallskip\noindent  
\noindent Chapter  VII.  Factoring polynomials. 

\noindent 1.  Factorial and\ separably factorial fields. 2.  Extensions of (separably) factorial fields. 3.  Seidenberg fields. 4.  The fundamental theorem of algebra. 5.  Notes. 

 \smallskip\noindent  
\noindent Chapter  VIII.  Commutative Noetherian rings.

\noindent 1.  The Hilbert basis theorem. 2.  Noether normalization and the Artin-Rees lemma. 3.  The Nullstellensatz. 4.  Tennenbaum's approach to the Hilbert basis theorem. 5.  Primary ideals. 6.  Localization. 7.  Primary decompositions. 8.~Lasker-Noether rings. 9.  Fully Lasker-Noether rings. 10.  The principal ideal theorem. 11.  Notes.  

 \smallskip\noindent  
\noindent Chapter  IX.  Finite dimensional algebras.

\noindent 1.  Representations. 2.  The density theorem. 3.  The radical and summands. 4.~Wedderburn's theorem, part one. 5.  Matrix rings and division algebras. \linebreak 6.~Notes. 

 \smallskip\noindent  
\noindent Chapter  X.  Free groups. 

\noindent 1.  Existence and uniqueness.  2.  Nielsen sets. 3.  Finitely generated subgroups of free groups.  4.  Detachable subgroups of finite-rank free groups.  5.  Conjugate subgroups. 6.  Notes.

 \smallskip\noindent  
\noindent Chapter  XI.  Abelian groups.

\noindent 1.  Finite-rank torsion-free groups. 2.  Divisible groups. 3.  Height functions on $p$-groups. 4.  Ulm's theorem. 5.  Construction of Ulm groups. 6.  Notes.

 \smallskip\noindent  
\noindent Chapter  XII.  Valuation theory.

\noindent 1.  Valuations.  2.  Locally precompact valuations. 3.  Pseudofactorial fields. 4.  Normed vector spaces. 5.  Real and complex fields. 6.  Hensel's lemma. 7.~Extensions of valuations. 8.  $e$ and $f$. 9.  Notes.

 \smallskip\noindent  
\noindent Chapter  XIII.  Dedekind domains.

\noindent 1.  Dedekind sets of valuations. 2.  Ideal theory. 3.  Finite extensions.

 \smallskip \noindent Bibliography. Index.

\medskip 
In the following sections we comment some significant examples
of classical \thos to which the \cov reformulation brings a new light and precise additional informations.

We also give some examples of  \thos which are trivial in \clama and yet very important from the \algq point of view.

\section{Principal ideal domains and \mtfs on these rings}

In \clama, a principal ideal domain is an integral ring in which all \ids are principal. From a constructive point of view, even the two-element field does not satisfy this \dfn: consider an ideal generated by a binary sequence; finding a generator of this ideal is the same thing as deciding if the sequence is identically zero, which amounts to LPO.  

An algorithmically relevant \dfn, classically \eqv to the classical one, is that of a discrete Bézout integral ring that satisfies a precisely formulated \noe condition.

\begin{cadre} \label{quote7}
A \textbf{GCD-monoid} is a cancellation [commutative] monoid in which each pair of
elements has a greatest common divisor. 
A \textbf{GCD-domain} is a discrete
domain whose nonzero elements form a GCD-monoid.

\dotfill\medskip 

A \textbf{principal ideal} of a commutative monoid $M$ is a subset $I$ of $M$
such that $I=Ma=\{ma:m\in M\}$ for some $a$ in $M$. We say that $M$
satisfies the \textbf{divisor chain condition} if for each ascending chain $I_{1}\subseteq
I_{2}\subseteq I_{3}\subseteq \cdots $ of principal ideals, there is $n$
such that $I_{n}=I_{n+1}$. 

A discrete domain is said to satisfy the divisor
chain condition if its monoid of nonzero elements does.

\smallskip \noindent \textbf{[IV.]2.7 Definition.} A \textbf{B\'{e}zout domain}%
\index{Bezout domain@B\'{e}zout domain} is a discrete domain such that for
each pair of elements $a,b$ there is a pair $s,t$ such that $sa+tb$ divides $%
a$ and $b$. A \textbf{principal ideal domain}%
\index{Principal ideal!domain} is a B\'{e}zout domain which satisfies the
divisor chain condition.
\end{cadre}


The classical structure \tho says that a \mtf on a PID
is a direct sum of a finite rank free submodule and of the torsion submodule, itself equal to a direct sum of modules  $R/(a_i)$ with the non-zero $a_i$ put in an order where each $a_i$ divides the next one.

The purest \algq form of this \tho is the \tho of reduction of a matrix into a Smith normal form.

\begin{cadre} \label{quote8}
A matrix $A=(a_{ij})$ is in \textbf{Smith normal form} if it is diagonal and $a_{ii}|a_{i+1,i+1}$ for
each $i.$

\medskip \noindent \textbf{[V.]1.2 Theorem.}
\emph{Each matrix over a principal ideal domain is equivalent
to a matrix in Smith normal form.}

\smallskip \noindent \textbf{[V.]1.4 Theorem.}
\emph{Two $m\times n$ matrices in Smith normal form over a
GCD-domain are equivalent if and only if corresponding elements are
associates.}
\end{cadre}

The structure \tho for \mpfs follows directly from \Tho V.1.2.

\begin{cadre} \label{quote9}
\noindent \textbf{[V.]2.3 Theorem} \textbf{(Structure \tho).} \emph{Let $M$ be a finitely presented module over a principal
ideal domain $R$. Then there exist principal ideals $I_{1}\supseteq
I_{2}\supseteq \cdots \supseteq I_{n}$ such that $M$ is isomorphic to the
direct sum $R/I_{1}\oplus R/I_{2}\oplus \cdots \oplus R/I_{n}$.}
\end{cadre}

Since the ring is discrete by \dfn, we can separate the sum into two pieces: the beginning, for indices from $1$ to $ k $ say, is the torsion submodule, with $I_k=(a_k)\neq0$, and the second piece, for $j>k$ with $a_j=0$, is a free module of rank $n-k$. On the other hand, in order to know which $I_j$'s ($j\leq k$) are equal to $R$ (and thus could be removed without damage), we need to have a test of invertibility for \elts of $R$, which in this case is equivalent to having a divisibility test between two \elts.

In \clama,  \Tho V.2.3 is stated for \mtfs. From a classical point of view the \mtfs over a PID are \pf, while from a \cof point of view  it is clearly impossible to have an \algo to achieve this implication, even in the simple case of the $\gZ$-module  $\gZ/I$ where $I$ is countably generated (e.g. generated by a binary sequence).

The way in which Bourbaki (Algebra, chapter VII) treats these \thos deserves to be compared.
The structure \tho is given before the Smith reduction \tho for matrices. And the \demo, which uses LEM, fails to produce an \algo to make the \tho explicit.

\section{Factorization \pbs}

Theorem IV.4.7 (i) below is usually shown for unique factorization domains, but the underlying Noetherian condition is in fact useless. 

\begin{cadre} \label{quote10}
\noindent \textbf{[IV.]4.7 Theorem.} {\it Let $R$ be a discrete domain.

 \noindent {\rm (i)}  If $R$ is a GCD-domain, then so is $R[X].$
}
\end{cadre}

The reader is invited to appreciate the elegance of the proof in [CCA].

\smallskip The classical \tho of factorization of an element into a product of prime factors in a GCD monoid  satisfying the divisor chain condition is inaccessible from an \algq point of view. It is replaced in \coma by a slightly more subtle \tho. This new \tho can generally be used instead of the classical one when needed to obtain concrete results.
 
\begin{cadre} \label{quote10a}
\noindent \textbf{[IV.]1.8 Theorem} \textbf{(Quasi-factorization).} \emph {Let $x_{1},\ldots ,x_{k}$ be elements of
a GCD-monoid $M$ satisfying the divisor chain condition. Then there is a
family $P$ of pairwise relatively prime elements of $M$ such that each $x_{i}
$ is an associate of a product of elements of $P.$}
\end{cadre}

\begin{cadre} \label{quote10b}
Let $M$ be a cancellation monoid. An element $a\in M
$ is said to be \textbf{bounded by $n$} if whenever $a=a_{0}\cdots a_{n}$
with $a_{i}\in M$, then $a_{i}$ is a unit for some $i$. An element of $M$ is 
\textbf{bounded} if it is bounded by $n$ for some $n\in \mathbf{N}$;
the monoid $M$ is \textbf{bounded} if each of its elements is bounded. A
discrete domain is \textbf{bounded} if its nonzero elements form a bounded
monoid.

A GCD-domain
satisfying the divisor chain condition is called a \textbf{quasi-UFD}.
\end{cadre}

The quasi-UFDs  and the bounded GCD-domains are two constructive versions (that are not constructively equivalent) of the classical notion of a UFD. In fact, we find in [CCA] still three other constructive versions of this classical notion.

\begin{cadre} \label{cadre}
\smallskip \noindent \textbf{[IV.]2.1 Definition.} A discrete domain $R$ is called a \textbf{unique factorization domain}, or \textbf{UFD}, if each nonzero element $r$ in $R$ is either a unit or has an
essentially unique factorization into irreducible elements, that is, if $r=p_{1}\cdots p_{m}$ and $r=q_{1}\cdots q_{n}$ are two factorizations of $r$
into irreducible elements, then $m=n$ and we can reindex so that $p_{i}\sim
q_{i}$ for each $i$. We say that $R$ is \textbf{factorial} if $R[X]$ is a UFD.

 \smallskip \noindent Call a
discrete field $k$ \textbf{fully factorial} if any finite-dimensional extension of $k$ is
factorial.

\end{cadre}

The five constructive versions are in classical mathematics equivalent to the classical notion, but they introduce algorithmically relevant distinctions, totally invisible in classical mathematics, due to the use of LEM, which annihilates these relevant distinctions.  In Theorem IV.4.7 the points (ii) (attached to the point (i)) and (vi) (i.e. (i) and (v)) are two distinct, inequivalent versions of the same classical theorem about UFDs.

\begin{cadre} \label{quote11}
\noindent \textbf{[IV.]4.7 Theorem.} {\it Let $R$ be a discrete domain.
\begin{enumerate}
\item[(i)]  If $R$ is a GCD-domain, then so is $R[X].$

\item[(ii)]  If $R$ is bounded, then so is $R[X].$

\item[(iii)]  If $R$ has recognizable units, then so does $R[X].$

\item[(iv)]  If $R$ has decidable divisibility, then so does $R[X].$

\item[(v)]  If $R$ satisfies the divisor chain condition, then so does $%
R[X]. $

\item[(vi)]  If $R$ is a quasi-UFD, then so is $R[X].$
\end{enumerate}
}
\end{cadre}

Concerning factorization problems for polynomials over a discrete field, the algorithmic situation is not correctly described by classical mathematics. E.g. factorization of polynomials in $k[X]$ where $k$ is a discrete field is not a trivial thing, contrarily to what is stated in \clama. 

Chapter VII of [CCA] explores the factorization problems in polynomial rings  in great detail.

The basic constructive theorem on this subject is given in Chapter VI. As it happens that the characteristic of a field or a ring is not known in advance, but can be revealed during a construction, some precautions are necessary in the statements, as below in  point (i). Note that if we discover a prime number~$p$ equal to zero in a ring $k$, it is necessarily unique (unless the ring is trivial).

In the following theorem, if $k$ is a discrete field, then we simply drop the alternative \gui{$k$ has a nonzero nonunit}. But it happens in [CCA] that the theorem is used in the precise form given here, e.g.\ in Chapter IX about the structure of finite-dimensional algebras.
\begin{cadre} \label{quote12}
\noindent \textbf{[VI.]6.3 Theorem.} {\it Let $k$ be a discrete commutative ring with
recognizable units, and $S$ a finite set of monic polynomials in $k[X]$.
Then either $k$ has a nonzero nonunit or we can construct a finite set $T$
of monic polynomials in $k[X]$ such that

\begin{enumerate}
\item[(i)]  Each element of $T$ is of the form $f(X^{q})$ where $f$ is
separable, and $q=1$ or $q$ is a power of a prime that is zero in $k$.

\item[(ii)]  Distinct elements of $T$ are strongly relatively prime.

\item[(iii)]  Every polynomial in $S$ is a product of polynomials in $T$.
\end{enumerate}
}
\end{cadre}

When $ k $ is a discrete field, we thus obtain, starting from a given family of univariate polynomials, a family of separable strongly relatively prime monic polynomials which gives a more precise version of the quasi-factorization theorem IV.1.8 (which deals with quasi-UFDs).

%

\section{\noe rings, primary decompositions and the principal ideal \tho}

An $R$-module is said to be \textbf{strongly discrete} if  \tf submodules are detachable.\footnote{In [CCA], the terminology is \gui{module with detachable submodules}, it was later replaced by \gui{strongly discrete module}. See e.g. \cite[Richman 1998]{Ric1998b}.} 
It is said to be \textbf{\coh}\footnote{Bourbaki (Algebra, Chapter X, or Commutative Algebra Chapter I) calls  \emph{pseudo coherent module} what [CCA] calls \coh module (as in quasi all texts in english literature), and \emph{\coh module} what [CCA] calls \pf \coh module. This is to be linked to \gui{Faisceaux Algébriques Cohérents} by \hbox{J.-P. Serre}.
Note also that the Stacks Project (Collective work,
\url{http://stacks.math.columbia.edu}) uses  Bourbaki's \dfn for \coh modules.} if any \tf submodule is \pf. The notion of \fdi coherent ring is fundamental from the algorithmic point of view in commutative algebra. In particular for the following reason: on a strongly discrete coherent ring, linear systems are perfectly understood and mastered.\footnote{In the article of Posur cited above, these rings are called \gui{computable}.}

In usual textbooks in \clama,  this notion is usually hidden behind that of a \emph{\noe} ring, and rarely put forward.
In \clama every \noe ring $R$ is \coh because every submodule of $R^{n}$ is \tf, and every \mtf is \coh for the same reason.
Furthermore, we have the Hilbert basis \tho, which states that \emph{if $R$ is \noe, then every \pf $R$-\alg is also a \noe ring,} whereas the same statement does not hold if one replaces \gui{\noe} with \gui{\coh} (see \cite[Soublin, 1970]{Sou70}).

From an \algq point of view however, it seems
impossible to find a satisfying \cov formulation 
of \noet which implies \cohc,
and \cohc is often the most important \prt from an \algq point of view.
Consequently, from a \cov point of view, \cohc must be added 
when we use the notion of a \noe ring or module.

The definition adopted for  \textbf{Noetherian module} in [CCA] is: a module in which any ascending chain of \tf submodules admits two equal consecutive terms. It is a constructively acceptable definition, equivalent in classical mathematics to the usual definition.

The classical \tho stating that over a \noe ring
every \tf \Amo is \noe is often advantageously replaced by the following \cof         \thos.

\noindent  \emph{Over a \coh  ring (resp. \fdi \coh)
every \pf \Amo is  \coh  (resp. \fdi \coh)}.

\noindent  \emph{Over a  \noe \coh ring
every \pf \Amo is   \noe \coh}.

\smallskip Two important classical results about Noetherian rings have constructive proofs within the framework given by [CCA].
\begin{cadre} \label{quote14}
\noindent \textbf{[VIII.]2.7 Theorem} \textbf{(Artin-Rees).} 
\emph{Let $I$ be a finitely generated ideal of a coherent
commutative Noetherian ring $R$. Let $N$ be a finitely generated submodule
of a finitely presented $R$-module $M$. Then there is $k$ such that for all $%
n\geq k$ we have\ 
\[
I^{n-k}(I^{k}M\cap N)=I^{n}M\cap N.
\]}
\noindent \textbf{[VIII.]2.8 Theorem} \textbf{(Krull intersection theorem).} 
\emph{Let $M$ be a finitely presented module
over a coherent commutative Noetherian ring $R$, and let $I$ be a finitely
generated ideal of $R$. Let $A=\cap _{n}I^{n}M$. Then $a\in Ia$ for each $%
a\in A$, so $IA=A$.}

\end{cadre}

\subsection*{Hilbert basis \tho}

Which are the \coh rings $R$ such that the \pol rings $R[\Xn]$ are also \coh? From a \cof point of view, we know two classes of rings sharing this \prt: \coh \noe rings (see below) and Prüfer domains  (see \cite[Yengui, 2015, Chapter 4]{Yen2015}).

The Hilbert basis \tho for the \dfn of \noet given in [CCA] is Theorem VIII.1.5 below. Proofs go back to 1974 (\cite[Richman, 1974]{Ric1974} and \cite[Seidenberg, 1974]{Sei1974}, see also \cite[Seidenberg, 1971]{Sei1971} and \cite[Seidenberg, 1973]{Sei1973} for \pol rings over a discrete field). These proofs are very clearly laid out in [CCA].

\begin{cadre} \label{quote15}
\noindent \textbf{[VIII.]1.5 Theorem  (Hilbert basis theorem).}
\emph{If $R$ is a coherent Noetherian ring, then so
is $R[X]$. If, in addition, $R$ has detachable left ideals,\footnote{i.e., is left \fdi} then so does 
$R[X]$.}
\end{cadre}

There is an analogous \tho in Computer Algebra  (see \cite[1994, \Tho 4.2.8]{AdLo1994}) saying that for a \coh \noe \fdi ring $R$, 
there is a  \gui{Gröbner basis \algo} computing the leading \id of a \itf in $R[\Xn]$ for a given monomial order.  In fact, this Computer Algebra \tho and \Tho VIII.1.5 are essentially the same result. One is easily deduced from the other.

Nevertheless we note that \algos for these theorems are quite different from each other. Moreover, authors in 1994 seem to ignore that the \pb was solved essentially in 1974, and \algos in \cite{AdLo1994} are not certified \cot (in fact, from the  \demo, no bound can be estimated for the number of steps as depending on the data.)   

\subsection*{Primary decomposition theorem}

 [CCA] gives an adequate \cof theory of primary decompositions.  
 This is based on the work of Seidenberg,  \cite[1978]{Sei1978} and \cite[1984]{Sei1984}. In [CCA] this work is made more simple and synthetic.

Let $R$ be a commutative ring. An ideal $Q$ of $R$ is said to be \textbf{primary} if $xy\in Q$ implies $x\in Q$ or $y^{n}\in Q$ for some 
$n$. One sees that $\sqrt{Q}$ is a \idep~$P$. 

[CCA] gives a variant w.r.t.\ the usual terminology, with no importance in the case of \noe rings for \clama:  \ids are all \tf.
A \textbf{primary decomposition} of an ideal $I$ in a commutative ring is a finite family of \tf primary \ids 
$Q_{1},\ldots ,Q_{n}$ such that the~$\sqrt{Q_i}$ are  \tf and $I=\bigcap_{i}Q_{i}$.
In this case the \id~$I$ is said to be \textbf{decomposable}.
  In \clama, every \id of a \noe ring has a primary decomposition.
  
In a \cof framewok, which convenient hypotheses do we have to add for a \coh \noe \fdi ring in order to get primary decompositions? 
A possible answer is the following one, given in [CCA]. 

\begin{cadre} \label{quote16}
A \textbf{Lasker-Noether} ring is a coherent Noetherian ring with
detachable ideals such that the radical of each finitely generated ideal is
the intersection of a finite number of finitely generated prime ideals.
\end{cadre}

This \dfn is \cot acceptable and applies to usual examples like $\gZ$,
$\gQ[X]$, and $k[X]$ when $k$ is an \agqt closed discrete field: they are  clearly \cot Lasker-Noether for this definition. Many other usual examples are also available, as explained below. 

In fact, when $k$ is a discrete field, $k[X]$ is easily seen to be Lasker-Noether \ssi  $k$ is a factorial field. This equivalence has no meaning  in \clama since all fields are factorial. Nevertheless it should be possible to state an analogous result for mechanical computations using Turing machines.

\smallskip The first \prts of Lasker-Noether rings are summarized in three theorems.

\begin{cadre} \label{quote17}
\noindent \textbf{[VIII.]8.1 Theorem.} \emph{Let $S$ be a multiplicative submonoid of a Lasker-Noether ring $R$ such that $I\cap S$ is either empty or nonempty for each finitely generated ideal $I$\ of $R$. Then $S^{-1}R$ is a
Lasker-Noether ring.}
\end{cadre}

If $S=R\setminus P$ for a \idep $P$, condition \gui{$I\cap S$ is either empty or nonempty} means that \gui{$I$ either is contained in $P$ or  is not}. Since $I$ is \tf, the test exists  \ssi $P$ is detachable. So, \tho VIII.8.1 implies that for each detachable \idep, and so for each  \tf \idep, the localization $R_P$ is Lasker-Noether.

\begin{cadre} \label{quote18}
\noindent \textbf{[VIII.]8.2 Theorem.} \emph{Let $R$ be a Lasker-Noether ring, and let $I$\ be a
finitely generated ideal of $R$. Then $R/I$ is a Lasker-Noether ring.}

\smallskip \noindent \textbf{[VIII.]8.5 Theorem  (Primary decomposition theorem).} \emph{Let $R$ be a Lasker-Noether ring. Then
each finitely generated ideal of $R$ has a primary decomposition.}
\end{cadre}

\subsection*{Principal ideal \tho}

A more elaborate \prt of Lasker-Noether rings is the famous principal ideal theorem of Krull and the fact that finitely generated proper prime ideals have a well-defined height. 
 
\begin{cadre} \label{quote19}
\noindent \textbf{[VIII.]10.4 Theorem   (Generalized principal ideal theorem).} \emph{Let $R$ be a Lasker-Noether ring. Let $%
I=(a_{1},\ldots ,a_{n})$. Then every minimal prime ideal over $I$\ has
height at most $n$.}

\smallskip \noindent \textbf{[VIII.]10.5 Theorem.} \emph{Let $P$ be a finitely generated proper prime ideal
of a Lasker-Noether ring $R$. Then there is $m$ such that $P$ has height $m$, and $P$ is a minimal prime over some ideal generated by $m$ elements.}
\end{cadre}

\subsection*{Fully Lasker-Noether rings}

Finally, it is important to give a \cov answer to the following: 
which convenient hypotheses do we have to add for a Lasker-Noether ring $R$ in order to get that 
$R[\Xn]$ is also Lasker-Noether?  

\begin{cadre} \label{quote20}
Call $R$ a \textbf{fully Lasker-Noether} ring if it is a Lasker-Noether ring and if for each finitely generated prime ideal $P$ of $R$, the field of quotients of $R/P$ is fully factorial. Note that the ring of integers 
$\mathbf{Z}$ is a fully Lasker-Noether ring, as is any fully factorial field.
\end{cadre}

The following three theorems (with the previous theorems about Lasker-Noether rings) show that in this context (i.e.\ with this constructively acceptable definition equivalent to the  definition of a Noetherian ring in classical mathematics), a very large number of classical theorems  concerning Noetherian rings now have a constructive proof and a clear meaning. It sounds like a \gui{miracle} of the same kind as Bishop's book.

\begin{cadre} \label{quote21}
\noindent \textbf{[VIII.]9.1 Theorem.} \emph{Let $I$\ be a finitely generated ideal of a fully Lasker-Noether ring $R$. Then $R/I$ is a fully Lasker-Noether ring.}

\smallskip \noindent \textbf{[VIII.]9.2 Theorem.} \emph{If $P$ is a detachable prime ideal of a fully Lasker-Noether ring $R$, then $R_{P}$ is a fully Lasker-Noether ring.}

\smallskip \noindent \textbf{[VIII.]9.6 Theorem.} \emph{If $R$ is a fully Lasker-Noether ring, then so is $R[X]$.}
\end{cadre}

\smallskip \noindent \emph{Note.} The paper \cite[Perdry, 2004]{Per2004} defines a notion of \noet which is \cot stronger than the one in [CCA]. The usual examples of \noe rings are  \noe in this meaning.
With this notion, the definition of a Lasker-Noether ring becomes more natural: it is a  \noe \coh \fdi ring in which we have a primality test for \itfs.
The paper gives a nice theory of fully Lasker-Noether rings in this  context.

\smallskip \noindent \emph{Note.} The computation of primary decompositions   in polynomial rings over discrete fields or over $\gZ$ is an active area of research in Computer Algebra. The seminal paper of Seidenberg is sometimes cited, but not the book [CCA].

\section{Wedderburn structure \tho   for finite-\hspace{0pt}dimensional  \klgs}

We deal here with unitary associative \klgs which are finite-dimensional \kevs on a \cdi $k$. In other words, these algebras are isomorphic to a \tf subalgebra of an algebra of matrices $ E_k (k ^ n) $ (the algebra of $k$-endomorphisms of the vector space $k^n$). We shorten the terminology by speaking of \gui{\klg of finite dimension}.

If  $A$ is a not necessarily commutative ring, its \textbf{Jacobson radical}  is the set $I$ of elements $x$ such that $1 + xA \subseteq  A^\times $. It is a  (two-sided) ideal and the Jacobson radical of the quotient $A / I$ is zero.

When $A$ is a $k$-algebra of finite dimension, this radical can also be defined as the \textbf{nilpotent radical}: $\rad (A)$ is the set of elements $x$ such that the left ideal $xA$ is nilpotent, i.e.\ there exists an integer $n$ such that every product $xa_1 \cdots xa_n$ is zero.

Let $A$ be a $k$-algebra of finite dimension. We can construct a basis of the center of $A$ as well as the minimal polynomial over $k$ of an arbitrary element of $A$. We can also construct a basis of the left ideal and another of the two-sided ideal generated by a finite part of $A$. But it may be difficult to construct a basis of the radical, and we cannot generally state that the radical is finite-dimensional (over $k$).

Nevertheless, we know how to construct objects whose counterparts are trivial in classical mathematics (if we do not try to construct them!). For example, as an alternative to the construction of the radical, we have the following theorem.
\begin{cadre} \label{quote22}
\noindent \textbf{[IX.]3.3 Theorem.} \emph{Let $A$\ be a finite-dimensional \klg and $L$ a
finite-dimensional (left) ideal of $A$. Then either $L\cap \rad A\neq 0$ or $%
A=L\oplus N$ for some (left) ideal $N$.}
\end{cadre}

A module $M$ is \textbf{reducible} if it has a nontrivial submodule---otherwise it is 
\textbf{irreducible} (or \textbf{simple}).

A \klg  is said to be \textbf{simple}  if each two-sided ideal is trivial. 
When the algebra is discrete (as in the present context)  the \dfn amounts to saying that if an element is nonzero, the (two-sided) ideal  it generates contains 1.

The first part of Wedderburn's structure theorem says that
every finite-dimensional $k$-algebra with zero radical is a product of simple algebras. Here is the \cov  reformulation given in [CCA]. A field $k$ is called separably factorial when separable \pols in $k[X]$ have a prime decomposition.  
 
\begin{cadre} \label{quote23}
We now characterize separably factorial fields in terms of decomposing
algebras into products of simple algebras. This is the first part of
Wedderburn's theorem.

\smallskip \noindent \textbf{[IX.]4.3 Theorem.} \emph{A discrete field $k$ is separably factorial if and only if every finite-dimensional $k$-algebra with zero radical is a product of simple algebras.}
\end{cadre}

A clarification concerning the ability to construct a basis of the radical is given in the following corollary.

\begin{cadre} \label{quote24}
\noindent \textbf{[IX.]4.5 Corollary.} \emph{A discrete field $k$ is fully factorial if and only if every finite-dimensional algebra $A$ over $k$ has a finite-dimensional nilpotent ideal $I$ such that $A/I$ is a product of simple $k$-algebras.}
\end{cadre}

The second part of Wedderburn's structure theorem for semi-simple algebras
says that a finite-dimensional simple algebra is isomorphic to a full ring
of matrices over a division algebra.

The \cov version of this \tho given in [CCA]  elucidates in a surprising way the computational content of this classical theorem.

\begin{cadre} \label{quote25}
\noindent \textbf{[IX.]5.1 Theorem.} 
 {\it Let $A$\ be a finite-dimensional $k$-algebra,
and $L$ a nontrivial left ideal of $A$. Then either
\begin{enumerate}
\item[(i)]  $A$\ has a nonzero radical

\item[(ii)]  $A$\ is a product of finite-dimensional $k$-algebras

\item[(iii)]  $A$\ is isomorphic to a full matrix ring over some $k$-algebra
of dimension less than $A$.
\end{enumerate}
}

\dotfill\medskip

The fundamental problem is to be able to recognize whether a given
finite-dimensional algebra is a division algebra or not, in the sense of
being able either to assert that it is a division algebra or to construct a
nontrivial left ideal. If we could do that, then Theorem  IX.5.1
would imply that every finite-dimensional $k$-algebra has a finite
dimensional radical, and modulo its radical it is a product of full matrix
rings over division algebras. This condition is equivalent to being able to
recognize whether an arbitrary finite-\hspace{0pt}dimensional representation of a
finite-dimensional $k$-algebra is reducible.

\smallskip \noindent \textbf{[IX.]5.2 Theorem.} {\it The following conditions on a discrete field $k$ are equivalent.
\begin{enumerate}
\item[(i)]  Each finite-dimensional $k$-algebra is either a division algebra
or has a nontrivial left ideal.

\item[(ii)]  Each finite-dimensional left module $M$ over a finite
dimensional $k$-algebra $A$ is either reducible or irreducible.

\item[(iii)]  Each finite-dimensional $k$-algebra $A$\ has a
finite-dimensional radical, and $A/\rad A$ is a product of full matrix rings
over division algebras.
\end{enumerate}
}
\end{cadre}

And we remain a little disappointed with these questions at the end of chapter IX.

\begin{cadre} \label{quote25a}
For what fields $k$ do the conditions of Theorem  5.2 hold?
Finite fields and algebraically closed fields provide trivial examples. The
field of algebraic real numbers admits only three finite-dimensional
division algebras, and a constructive proof of this statement shows that
this field satisfies the conditions of Theorem 5.2.

\smallskip \noindent \textbf{[IX.]5.3 Theorem.} \emph{Let $k$ be a discrete subfield of $\gR$ that is algebraically closed in $\gR$, and $H=k(i,j)$ the quaternion algebra over $k$.
If $A$ is a finite-dimensional algebra over $k$, then either $A$\ has a
zero-divisor, or $A$\ is isomorphic to $k$, to $k(i)$, or to $H$.}

\smallskip 
Does the field $\mathbf{Q}$ of rational numbers satisfy the conditions of
Theorem  5.2? Certainly we are not going to produce a
Brouwerian counterexample when $k=\mathbf{Q}$. Probably a close analysis of
the classical theory of division algebras over $\mathbf{Q}$, in analogy with
Theorem  5.3, will yield a proof.\medskip
\end{cadre}

\section{Dedekind domains}

\begin{cadre} \label{quote26a}
Although it is commonly felt that algebraic number theory is essentially
constructive in its classical form, even those authors who pay particular
attention to the constructive aspects of the theory employ highly
nonconstructive techniques which nullify their efforts. In
\cite[Borevich-Shafarevich, 1966]{BS66}, for example, it is assumed that every
polynomial can be factored into a product of irreducible polynomials (every
field is factorial) and that given a nonempty subset of the positive
integers you can find its least element.
\end{cadre}

The constructive theory of Dedekind domains in [CCA] allows us to give an explicit version of the classical statements of number theory and algebraic geometry concerning local fields, for example in the book of
J.-P. Serre \cite{CorpsLocaux}. This theory also gives the appropriate hypotheses to account for the classical results concerning Dedekind domains, as found, for example, in Bourbaki.

This requires giving sufficiently precise and binding definitions, beginning with those in the theory of (rank-one) valuations.

For example, let us see the definitions concerning  Dedekind domains.

\begin{cadre} \label{quote26}
\noindent \textbf{[XIII.]1.1 Definition.} A nonempty discrete set $S$\ of nontrivial
discrete valuations on a Heyting field $k$ is a \textbf{Dedekind set}%
\index{Dedekind!set of valuations} if
\begin{enumerate}
\item[(i)]  For each $x\in k$ there is a finite subset $T$ of $S$ so that $%
|x|_{p}\leq 1$ for each $p\in S\setminus T$.

\item[(ii)]  If $q$ and $q^{\prime }$ are distinct valuations of $S$, and $%
\varepsilon >0$, then there exists $x\in k$ with $|x|_{p}\leq 1$ for each $%
p\in S$, such that $|x-1|_{q}<\varepsilon $ and $|x|_{q^{\prime
}}<\varepsilon $. Hence distinct valuations are inequivalent.
\end{enumerate}

\medskip Let $S$ be a Dedekind set of valuations on a Heyting field $k$. If $p\in S$,
then, because $p$ is nonarchimedean, the set $R(p)=\{x\in k:|x|_{p}\leq 1\}$
is a ring, which is local as $p$ is discrete. We call $R(p)$ the \textbf{%
local ring at $p$}. The elements of the ring $\bigcap _{p\in S}R(p)$ are
called the \textbf{integers at $S$}. A ring is a \textbf{Dedekind domain}%
\index{Dedekind!domain} if it is the ring of integers at a Dedekind set of
valuations on a Heyting field.
\end{cadre}

If the strong point is to give a constructive account of most of the classical theorems, a weak point is that for example a PID is a   Dedekind domain only in the case where we have algorithms of factorization of \idps into prime ideals.

We can compare this for example with the exposition in \cite{CACM}, where a definition is given that is constructively weaker but closer to the usual classical definition (see Definition XII-7.7 and Theorem XII-7.9). In \cite{CACM}, Dedekind domains have quasi-factorization of finite sets of \itfs, and  the total factorization Dedekind domains  correspond to the discrete Dedekind domains of [CCA].

\subsection*{Acknowledgement}
I am indebted to Thierry Coquand and Stefan Neuwirth for many relevant comments and suggestions.

\newpage 
\addcontentsline{toc}{section}{Références}
\small
\bibliographystyle{plain}
\bibliography{BibACMC}

\end{document}